\documentstyle[12pt]{article}
\begin{document}
\newtheorem{thm}{Theorem}[section]
\newtheorem{lem}{Lemma}[section]
\newtheorem{rem}{Remark}[section]
\newtheorem{prop}{Proposition}[section]
\newtheorem{cor}{Corollary}[section]
\newtheorem{question}{Question}
\title
{Remark on some conformally invariant integral equations:
 the method of moving spheres
}
\author{ YanYan Li\\ Department of Mathematics\\ Rutgers University\\
110 Frelinghuysen Rd.\\
Piscataway, NJ 08854
}
\date{}
\maketitle
\newcommand {\Bbb}{}

\setcounter{section}{0}

\section{Introduction}

For $n\ge 3$, consider
\begin{equation}
-\Delta u=n(n-2) u^{ \frac{n+2}{n-2}},
\qquad \mbox{on} \quad \Bbb {R}^n.
\label{eq1new}
\end{equation}
It was proved by Gidas, Ni and Nirenberg
\cite{GNN1} that any positive $C^2$ solution of
(\ref{eq1new}) satisfying
\begin{equation}
\liminf_{|x|\to \infty}\left( |x|^{n-2}u(x) \right)<\infty,
\label{finite}
\end{equation}
must be of the form
$$
u(x)\equiv \left( \frac a{ 1+a^2|x-\bar x|^2 }\right)^{ \frac {n-2}{2} },
$$
where $a>0$ is some constant and $\bar x\in \Bbb R^n$.

Hypothesis (\ref{finite})
was removed by Caffarelli, Gidas and Spruck in \cite{CGS};
 this is important
for applications.
Such Liouville type theorems have  been extended to general
conformally invariant fully nonlinear equations by
Li and Li (\cite{LL1}-\cite{LL4});
 see also related works of Viaclovsky (\cite{V3}-\cite{V2})
 and Chang, Gursky and Yang (\cite{CGY2}-\cite{CGY3}).
The method used in \cite{GNN1}, as well as 
in much of the above cited work, 
 is the method of moving planes.
The method of moving planes has become a very powerful tool
in the study of nonlinear elliptic equations,
see Alexandrov \cite{Al}, Serrin \cite{Serrin}, Gidas, Ni and Nirenberg
 \cite{GNN1}-\cite{GNN2},
Berestycki and Nirenberg \cite{BN}, and others.

In \cite{LZhu}, Li and Zhu
gave a proof of the above mentioned theorem of
Caffarelli, Gidas and Spruck
using the method of moving spheres
(i.e. the method of moving planes together with 
the conformal invariance), which fully exploits
the conformal invariance of the problem and,
as a result,  captures the solutions
directly rather than going through
the
usual procedure of proving radial symmetry of solutions  and then classifying
radial solutions. 
 Significant simplifications to the proof in
 \cite{LZhu} have been made in Li and Zhang \cite{LZhang}.
The method of moving spheres has been  used in \cite{LL1}-\cite{LL4}.

Liouville type theorems for various conformally
invariant equations have received much attention, see,
in addition to the above cited papers,
\cite{GS}, \cite{CL2}, \cite{CY1},
\cite{lin1}, \cite{WX}
and \cite{X}.

In this paper we study some conformally invariant integral equations.
Lieb proved in \cite{Lieb}, among other things, 
that there exist maximizing functions, $f$,
for the Hardy-Littlewood-Sobolev inequality
on $\Bbb R^n$:
$$
\| \int_{\Bbb R^n} \frac {f(y)}
{ |\cdot-y|^\lambda }dy\|_{ L^q(\Bbb R^n) }
\le N_{p,\lambda, n}\|f\|_{ L^p(\Bbb R^n) },
$$
with $N_{p,\lambda, n}$ being the sharp constant
and 
$\frac 1p+\frac \lambda n=1+\frac 1q$,
$1<p, q, \frac n\lambda<\infty$, $n\ge 1$.
When $p=q'=\frac q{q-1}$ or $p=2$ or $q=2$,
$N_{p,\lambda, n}$ and the maximizing $f'$s are 
explicitly evaluated.  When $p=q'$,
i.e., $p=\frac{2n}{2n-\lambda}$ and $q=\frac {2n}\lambda$,  
the Euler-Lagrange equation for a maximizing $f$ is,
modulo a positive constant multiple,
\begin{equation}
f(x)^{p-1}=\int_{\Bbb R^n}
\frac{ f(y) }{ |x-y|^\lambda }dy.
\label{lieb1}
\end{equation}  
Writing $\lambda=n-\alpha$ and $u=f^{p-1}$, then $0<\alpha
<n$,  and equation (\ref{lieb1}) becomes 
\begin{equation}
u(x)=\int_{\Bbb R^n} \frac{ u(y)^{\frac {n+\alpha}{n-\alpha} }  }
{  |x-y|^{n-\alpha}  }dy, \quad
\forall\ x\in\Bbb R^n.
\label{1aa}
\end{equation}

As mentioned above, maximizing solutions $f$ of 
(\ref{lieb1}) are classified in \cite{Lieb} and they are,
in terms of $u$, of the form  
\begin{equation}
u(x)\equiv  (\frac a {d+|x-\bar x|^2})^{ \frac {n-\alpha}2 },
\label{formnew}
\end{equation}
where $a, d>0$ and
some  $\bar x\in \Bbb R^n$.  Of course, $a$ is a fixed constant depending
only on $n$ and $\alpha$, while $d$ and $\bar x$ are
free.

Equation (\ref{1aa}), or (\ref{lieb1}),
is conformally invariant in the following  sense.
Let $v$ be a positive function on $\Bbb R^n$, for
  $x\in \Bbb R^n$ and $\lambda>0$,
we define
\begin{equation}
v_{x,\lambda}(\xi)=(\frac \lambda{|\xi-x|})^{ n-\alpha}
v(\xi^{x,\lambda}),\qquad\xi\in \Bbb R^n,
\label{v24}
\end{equation}
where
$$
\xi^{x,\lambda}=x+ \frac {\lambda^2(\xi-x)}{|\xi-x|^2}.
$$
Then, if $u$ is a solution of (\ref{1aa}), 
so is $u_{x,\lambda}$ for any $x\in \Bbb R^n$ and $\lambda>0$.
The conformal  invariance of (\ref{1aa}) was used 
in \cite{Lieb}.  More studies on issues
concerning the Hardy-Littlewood-Sobolev
inequality, among
other things,  were given by Carlen and Loss in \cite{CL}-\cite{CL5},
where the conformal invariance of the problem was further exploited.

After classifying all maximizing
solutions of (\ref{lieb1}),
Lieb raised the beautiful question (page 361 of \cite{Lieb})
on the (essentially) uniqueness of solutions of (\ref{lieb1}),
or,  equivalently, of (\ref{1aa}).
He produced (page 363 of \cite{Lieb}) a nontrivial
$2n$ parameter family of solutions of 
equation (\ref{lieb1}), or 
(\ref{1aa}), has other solutions
which are not as regular as the maximizers.
For instance, modulo a positive constant,
$ |x|^{ \frac {\alpha-n}2 }$ is a solution
of of (\ref{1aa}).

In a recent paper, Chen, Li and Ou 
established the following  result which answers the
question of Lieb in the class of $L^\infty_{loc}(\Bbb R^n)$.
\begin{thm} (\cite{CLO})  
Let  $u\in L^\infty_{loc}(\Bbb R^n)$ be
a positive  function
satisfying
(\ref{1aa}).
Then $u$ is given by
(\ref{formnew}) for some constants $a, d>0$ and 
some  $\bar x\in \Bbb R^n$.
\label{thm1new}
\end{thm}

In an earlier version of the present paper \cite{Li}, 
we gave  a simpler proof of Theorem \ref{thm1new}.
The proof, in the spirit of \cite{LZhu} and \cite{LZhang}
and following  Section 2 of \cite{LZhang},
fully exploits the conformal invariance
of the integral equation.
It is different from the one in \cite{CLO}.
In particular,
we do not follow the usual 
procedure of proving radial symmetry of solutions and then 
classifying radial solutions, and we do not need
to distinguish $n\ge 2$ and $n=1$.
This proof  is presented in Section 2.

Lieb pointed out to us that his question also concerns
functions which are not in $L^\infty_{loc}(\Bbb R^n)$.
In particular, it is not known a priori
that maximizers are in $L^\infty_{loc}(\Bbb R^n)$.
This has led us to study the question further and
to establish

\begin{thm} 
For $n\ge 1$, $0<\alpha<n$, let $u\in
L^{   
\frac{2n}{n-\alpha}
}_{loc}(\Bbb R^n)$ be a positive solution of 
(\ref{1aa}).  Then $u\in C^\infty(\Bbb R^n)$.
\label{thm1aa}
\end{thm}

An answer to the question of Lieb is therefore known
in the class of $L^{ \frac{2n}{n-\alpha} }_{loc}(\Bbb R^n)$.
The above mentioned solution
$ |x|^{ \frac {\alpha-n}2 }$ does not belong to
$L^{ \frac{2n}{n-\alpha} }_{loc}(\Bbb R^n)$, though it belongs
to $L^t_{loc}(\Bbb R^n)$ for any $t< \frac{2n}{n-\alpha}$.
The question remains unanswered
for the class of $L^t_{loc}(\Bbb R^n)$ for
$t<\frac{2n}{n-\alpha}$.

In the process of
proving Theorem \ref{thm1aa}, we have established the following
result 
 which should be of independent interest.

For $n\ge 1$ and $0<\alpha<n$, let  $V\in L^{\frac n\alpha}(B_3)$ 
be a non-negative function, set
\begin{equation}
\delta(V):=\|V\|_{ L^{\frac n\alpha}(B_3) }.
\label{P}
\end{equation}
\begin{thm} 
For $n\ge 1$, $0<\alpha<n$, $\nu> r>\frac n{n-\alpha}$, there exist
positive constants 
$\bar\delta<1$ and $C\ge 1$, depending only on $n, \alpha, r$ and $\nu$,
such that for any   $0\le V\in L^{\frac n\alpha}(B_3)$, 
with $\delta(V)\le \bar\delta$,
  $h\in L^\nu(B_2)$ and $0\le u\in L^r(B_3)$ satisfying
\begin{equation}
u(x)\le \int_{ B_3} \frac {V(y) u(y) }{  |x-y|^{ n-\alpha} }dy+h(x),
\qquad x\in B_2,
\label{linear}
\end{equation}
we have
\begin{equation}
\|u\|_{ L^\nu(B_{\frac 12}) } 
\le C\left(\|u\|_{ L^r(B_3) }+\|h\|_{ L^\nu(B_2) }\right).
\label{estimate}
\end{equation}
\label{thm5}
\end{thm}

\begin{cor} For $n\ge 1$, $0<\alpha<n$, $\nu> r>\frac n{n-\alpha}$,
$R_2>R_1>0$, let
 $0\le V\in L^{\frac n\alpha}(B_{R_2})$,
$h\in L^\nu(B_{R_1})$ that $0\le u\in L^r(B_{R_2})$ satisfy
$$
u(x)\le\int_{ B_{R_2} }  \frac {V(y) u(y) }{  |x-y|^{ n-\alpha} }dy+h(x),
\qquad x\in B_{R_1}.
$$
Then, for some $\epsilon>0$,
$u\in L^\nu(B_\epsilon)$.
\label{cor1}
\end{cor}

For $\alpha=2$ and $n\ge 3$, Theorem 
 \ref{thm5} is essentially equivalent to a result of Brezis and Kato
(Theorem 2.3 in \cite{BK}), so it  
 can  be viewed as an 
integral equation analogue of their theorem.
After informing Brezis of Theorem
 \ref{thm5}, he kindly pointed out that
it is similar to, though not the same as,
Lemma A.1 in \cite{BL}.  Indeed, our proof of
the theorem makes use of the explicit form of
the potential $|x|^{\alpha-n}$, and it is not clear to us
at this point whether the conclusion of 
the theorem still holds when
replacing $|x|^{\alpha-n}$ by any $Y\in L^{ \frac n{n-\alpha} }_w$,
the weak $ L^{ \frac n{n-\alpha} }$ space, as in 
Lemma A.1 of \cite{BL}.
Theorem \ref{thm1aa}, Theorem \ref{thm5} and Corollary \ref{cor1}
are established in Section 2.

We also study some equations similar to
(\ref{1aa}), though they do not have the same kind of conformal
invariance property.
For $n\ge 1$, $0<\alpha<n$ and $\mu>0$, 
let $u$ be positive Lebesgue measurable function in
$\Bbb R^n$
satisfying
\begin{equation}
u(x)=\int_{\Bbb R^n} \frac{ u(y)^\mu  }
{  |x-y|^{n-\alpha}  }dy, \quad
\forall\ x\in\Bbb R^n.
\label{1a}
\end{equation}

\begin{thm} Let  $n\ge 1$ and $0<\alpha<n$.
Then

(i) For $0<\mu< \frac n{n-\alpha}$,
equation (\ref{1a}) does not have any positive
Lebesgue measurable solution $u$, unless $u\equiv \infty$.

(ii) For  $ \frac n{n-\alpha}\le  \mu<\frac  {n+\alpha}{n-\alpha}$,
equation 
 (\ref{1a}) does not have any positive
solution $u\in  L^{  \frac {n (\mu-1) }\alpha  }_{loc}(\Bbb R^n)$.
\label{thm1a}
\end{thm}

For $\mu>\frac{n+\alpha}{n-\alpha}$,
we do not know whether  (\ref{1a})
has any positive solutions.
We know from Lemma \ref{lem0a} that if $u$ is a positive solution
in $L^{  \frac {n (\mu-1) }\alpha  }_{loc}(\Bbb R^n)$,
$u$ must be in $C^\infty(\Bbb R^n)$.
Theorem \ref{thm1a} is proved in Section 4.

In \cite{LL1}-\cite{LL4}, all conformally invariant second order fully nonlinear
equations are classified and Liouville type theorems 
are established for the elliptic ones.
It would be  interesting to identify as many as possible conformally invariant
integral equations for which  (essentially) uniqueness of solutions
can be obtained.
One class of such equations, similar to 
(\ref{1aa}), is
$$
u(x)=\int_{\Bbb R^n}
|x-y|^p u(y)^{-\frac {2n+p}p }dy, \quad
\forall\ x\in\Bbb R^n,
$$
where
$n\ge 1$ and $p>0$.
We study more general equations, similar to (\ref{1a}),
 including those which are not conformally invariant.
  
For $n\ge 1$, $p, q>0$, let $u$ 
be a non-negative Lebesgue measurable function in $\Bbb R^n$ satisfying
\begin{equation}
u(x)=\int_{\Bbb R^n}
|x-y|^p u(y)^{-q}dy, \quad
\forall\ x\in\Bbb R^n.
\label{1}
\end{equation}
\begin{thm}  For $n\ge 1$, $p>0$ and $0<q\le 1+\frac {2n}p$,
let  $u$ 
be a non-negative Lebesgue measurable function in $\Bbb R^n$ satisfying
(\ref{1}).
Then $q=1+\frac {2n}p$ and,
for some constants $a, d>0$ and 
some  $\bar x\in \Bbb R^n$,
\begin{equation}
u(x)\equiv \left( \frac {d+|x-\bar x|^2}a\right)^{ \frac p2}.
\label{form}
\end{equation}
\label{thm1}
\end{thm}

The proof of Theorem \ref{thm1}, similar to our proof 
of  Theorem \ref{thm1new}, is given in Section 5.
It turns out that for $n=3$, $p=1$ and $q=7$,
integral equation ({1})
is associated with some fourth order
conformal covariant operator on  $3-$dimensional
compact Riemannian manifolds, arising from 
 the study of conformal geometry. See, e.g., 
Paneitz \cite{P},
 Fefferman and Graham \cite{FG}, Branson \cite{Br-1}
 and Chang and Yang \cite{CY0}.
\begin{question}
Is equation ({1}), in the case $n\ne 3$,
$p>0$  and $q=1+\frac {2n}p$,
  is associated with some 
kind of pseudo-differential conformal covariant operators on $n-$dimensional
compact Riemannian manifolds, the same way the case
$n=3$, $p=1$ and $q=7$ is associated with the above mentioned
 fourth order
conformal covariant operator?
\end{question}

After posting \cite{Li} on the Archive and essentially
completing the proof of Theorem \ref{thm1}, we became aware
of some recent work of Xu \cite{X2} where he  proved Theorem \ref{thm1}
in the special case $n=3$, $p=1$
and $u\in C^4(\Bbb R^3)$.  He also proved in the same
paper that for $n=3$, $p=1$
and $q>7 (=1+\frac {2n}p)$,
 equation (\ref{1})
does not admit any non-negative solution
$u$ in $C^4(\Bbb R^3)$. 

\begin{question}
Is it true that for all $n\ge 1$, $p>0$ and $q>1+\frac {2n}p$
equation (\ref{1}) does not admit any positive solutions?
\end{question}

\section{ Proof of Theorem \ref{thm5}, Corollary \ref{cor1} and
 Theorem \ref{thm1aa} }
In this section we prove Theorem \ref{thm5}.
Let
$$
\xi(x):=\int_{ B_3} \frac {V(y) u(y) }{  |x-y|^{ n-\alpha} }dy
+h(x)-u(x)\ge 0,\qquad x\in B_2.
$$
Then
\begin{equation}
u(x)=(Lu)(x)
+
f(x)+h(x)-\xi(x),
\qquad x\in B_2,
\label{29-0new}
\end{equation}
where
$$
(Lu)(x)=\int_{ B_2} \frac {V(y) u(y) }{  |x-y|^{ n-\alpha} }dy,
\qquad x\in B_2,
$$
and
$$
f(x)=\int_{ 2<|y|<3 } \frac {V(y) u(y) }{  |x-y|^{ n-\alpha} }dy.
$$

Let $p$ be determined by $\frac 1r=\frac 1p-\frac \alpha n$,
then $p>1$ 
 and therefore,  by the property of the Riesz potential
(see, e.g., Theorem 1 on page 119 of \cite{Stein}),
\begin{eqnarray}
\|Lu\|_{ L^r(B_2)}
&\le
&
 C\| V u\|_{ L^p(B_2)}
=C\left( \|V^p u^p\|_{ L^1(B_2)}\right)^{ \frac 1p}
\nonumber\\
&\le & C \left( \|V^p \|_{
L^{ \frac r{r-p} }(B_2) }
\|u^p\|_{ L^{\frac rp}(B_2)  }  \right)^{ \frac 1p}
\le C \|V\|_{ L^{\frac n\alpha}(B_2)  }  \|u\|_{ L^r(B_2) },
\label{Lnew}
\end{eqnarray}
where $C$ depends on $\alpha, n$ and $r$.
Similarly
\begin{equation}
\|f\|_{ L^r(B_2) }
\le C \|V\|_{ L^{\frac n\alpha}(B_3)  }  \|u\|_{ L^r(B_3) }.
\label{fnew}
\end{equation}
It follows, using also the fact $u, \xi\ge 0$, that
\begin{equation}
\|\xi\|_{ L^r(B_2) }
\le C \|V\|_{ L^{\frac n\alpha}(B_3)  }  \|u\|_{ L^r(B_3) }
+C\|h\|_{ L^r(B_2) }.
\label{xinew}
\end{equation}

For $i=1,2,\cdots$, let
$$
G_i(z)=\min\left( \frac 1{|z|^{n-\alpha}}, i\right),
\qquad u_i(z)=\min\left( u(z), i\right),
$$
$$
\xi_i(x)=\min\left( \xi(x), i\right),
\qquad\mbox{and}\qquad f_i(x)=\int_{ 2<|y|<3 } G_i(x-y) V(y) u(y)dy.
$$

Some preliminary estimates on $\{f_i\}$:
\begin{lem} There exists some constant $C$, depending only on 
$n$ and $\alpha$,
such that
\begin{equation}
\|f_i\|_{ L^\infty(B_1) }
\le C\|u\|_{ L^r(B_3) }, \qquad  \|f_i\|_{ L^r(B_2) }
\le C\|u\|_{ L^r(B_3) }.
\label{indepedentnew}
\end{equation}
Moreover, for any $p<r$,
\begin{equation}
\lim_{i\to\infty} \|f_i-f\|_{ L^p(B_2) }=0.
\label{convergenew}
\end{equation}
\label{lemfnew}
\end{lem}

\noindent{\bf Proof of Lemma \ref{lemfnew}.}\
The first inequality in (\ref{indepedentnew})  follows easily:
$$
\|f_i\|_{ L^\infty(B_1) }
\le \|f\|_{ L^\infty(B_1) }
\le  C(n,\alpha) \int_{ 2<|y|<3} V(y)u(y)dy
\le  C(n,\alpha) \|u\|_{ L^r(B_3) }.
$$
Note that we have used the hypothesis  $\|V\|_{ L^{\frac n\alpha}(B_3) }
\le \bar \delta<1$.
The second inequality in (\ref{indepedentnew}) 
follows from (\ref{fnew}).

By the Fubini theorem,
$$
\lim_{i\to \infty}
\|f_i-f\|_{ L^1(B_2) }
\le \lim_{i\to \infty}
\|G_i(\cdot)-\frac 1{ |\cdot|^{n-\alpha} }\|_{ L^1(B_5) }
\int_{2<|y|<3}
V(y)u(y)dy=0.
$$
We deduce (\ref{convergenew})
from this and the second inequality in (\ref{indepedentnew})
using H\"older inequality.

\vskip 5pt
\hfill $\Box$
\vskip 5pt

Consider  the following integral equation of $w$,
\begin{equation}
w(x)=(L_iw)(x)
+f_i(x)+h(x)-\xi_i(x),
\qquad x\in B_2,
\label{29-1new}
\end{equation}
where
$$
(L_iw)(x):=
\int_{ |y|<2} 
G_i(x-y) V(y) w(y)dy.
$$

\begin{lem}  For $r\le q\le \nu$,  there exist some $0<\bar\delta<1$
and $C\ge 1$, 
depending only on $\alpha, n, r$ and $q$,  such that
if
 $0<\delta(V)\le \bar\delta$, then, for all $i$,  there exists
$w_i\in L^q(B_2)$ solving 
(\ref{29-1new}) with $w=w_i$,
satisfying
\begin{equation}
\|w_i\|_{ L^r(B_2) }
 +\|w_i^+\|_{ L^q(B_{\frac 12}) }\le C(\|u\|_{ L^r(B_3) }
+\|h\|_{ L^r(B_2) }),
\label{qenew}
\end{equation}
where $w_i^+(x)=\max( w_i(x), 0)$.
\label{lem4-2new}
\end{lem}

\noindent{\bf Proof of Lemma \ref{lem4-2new}.}\
Define, for $w\in L^q(B_2)$,
$$
(T_iw)(x)= (L_iw)(x)+f_i(x)+h(x)-\xi_i(x),
\qquad x\in B_2.
$$
Clearly, $L_i, T_i: L^q(B_2)\to L^q(B_2)$.

Let $p$ be determined by $\frac 1q=\frac 1p-\frac \alpha n$,
then, using the property of the Riesz potential as in (\ref{Lnew}),
$$
\|L_i w\|_{ L^q(B_2)}
\le
\|L(|w|)\|_{ L^q(B_2) }\le
 C \|V\|_{ L^{\frac n\alpha}(B_2)  }  \|w\|_{ L^q(B_2) }
\le C\bar\delta \|w\|_{ L^q(B_2) }.
$$
Here and below  (various) constant $C\ge 1 $
depends only on $r, q, \alpha$ and $n$.
Thus
\begin{equation}
\|T_i w\|_{ L^q(B_2) }
\le C\bar\delta \|w\|_{ L^q(B_2) }+\|f_i\|_{ L^q(B_2)}
+\|h\|_{ L^q(B_2)}+\|\xi_i\|_{ L^q(B_2)},
\label{thusnew}
\end{equation}
and
$$
\|T_i(w-v)\|_{ L^q(B_2) }
\le \|L_i(w-v)\|_{ L^q(B_2) }\le
C\bar\delta \|w-v\|_{ L^q(B_2) }.
$$

Fix some positive $\bar\delta$ with $ C\bar\delta
\le \frac 12$ and
set
$$
E_i=\{ w\in L^q(B_2)\ |\
\|w\|_{ L^q(B_2) }\le 2(\|f_i\|_{ L^q(B_2) } +
\|h\|_{ L^q(B_2) }+\|\xi_i\|_{ L^q(B_2)}) \}\subset L^q(B_2).
$$
Then, $T_i$ maps $E_i$ to itself and is a
contraction map.  So there exists some $w_i\in E_i$ such that
$T_i (w_i)=w_i$, i.e.,
\begin{equation}
w_i(x)=
\int_{ |y|<2} G_i(x-y) V(y)w_i(y)dy+f_i(x)+h(x)-\xi_i(x),
\qquad x\in B_2.
\label{qqnew}
\end{equation}

Taking
$q=r$ in (\ref{thusnew}), we
 obtain from
(\ref{qqnew}) and  (\ref{indepedentnew})
that
$$
\|w_i\|_{ L^{  r }(B_2) }
\le  \frac 12  \|w_i\|_{ L^r(B_2) }+\|f_i\|_{ L^r(B_2)}
+\|h\|_{ L^r(B_2)}+\|\xi\|_{ L^r(B_2)}.
$$
The estimate of
$\|w_i\|_{ L^r(B_2) }$
in (\ref{qenew}) follows from this,
in view of (\ref{xinew}) and the second inequality in (\ref{indepedentnew}).

Next we establish the second inequality in (\ref{qenew}).
For $0<t<s<1$, we have, by (\ref{qqnew}),
$$
w_i^+(x)\le I_i(x)+II_i(x)+f_i(x)+h(x),
$$
where
$$
I_i(x)=\int_{ |y|<s }
\frac{  V(y)w_i^+(y)  }
{  |x-y|^{n-\alpha}  }dy,
$$
and
$$
II_i(x)=\int_{s<|y|<2} \frac{  V(y)w_i^+(y)  }
{  |x-y|^{n-\alpha}  }dy.
$$

 By the property  of the Riesz potential,  
\begin{eqnarray*}
\|I_i\|_{ L^q(B_t) }
&\le &
C\|Vw_i^+\|_{ L^p(B_s) }
\le C \|V\|_{  L^{ \frac n\alpha }(B_s) } \|w_i^+\|_{ L^q(B_s) }
\\
&\le& 
  C\bar\delta \|w_i^+\|_{ L^q(B_s) }\le \frac 12  \|w_i^+\|_{ L^q(B_s) }.
\end{eqnarray*}
Using the estimate of $\|w_i\|_{ L^r(B_2) }$ in (\ref{qenew}),
\begin{eqnarray*}
\|II_i\|_{ L^q(B_t) }
&\le &
 C  (s-t)^{\alpha-n}  \int_{ s<|y|<2}
V(y)w_i^+(y) dy\\
&\le &  C(s-t)^{\alpha-n}
\|w_i\|_{ L^r(B_2) } \le
 C  (s-t)^{\alpha-n}(\|u\|_{ L^r(B_3) }+
\|h\|_{ L^r(B_2) }).
\end{eqnarray*}

With (\ref{indepedentnew})  and the above estimates,
we have, for all $0<t<s<1$,
$$
\|w_i^+\|_{ L^q(B_t) }\le \frac 12\|w_i^+\|_{ L^q(B_s) }
+ C (s-t)^{\alpha-n}(\|u\|_{ L^r(B_3) }+
\|h\|_{ L^r(B_2) }).
$$
By a calculus lemma (see, e.g., page 32 of \cite{GG}),
we have, for a possibly larger $C$,
still depending only on $r, q, \alpha$ and $n$, 
$$
\|w_i^+\|_{ L^q(B_t) }\le C (s-t)^{\alpha-n}(\|u\|_{ L^r(B_3) }+
\|h\|_{ L^r(B_2) }),
\qquad \forall\ 0<t<s<1.
$$
The estimate of $\|w_i^+\|_{ L^q(B_{\frac 12}
) }$ in   (\ref{qenew}) follows from the above.
Lemma \ref{lem4-2new} is established.

\vskip 5pt
\hfill $\Box$
\vskip 5pt

\noindent{\bf Proof of Theorem \ref{thm5}.}\
For any $r<q\le \nu$,
let $\bar\delta>0$ and  $\{w_i\}\in L^q(B_2)$ be given by Lemma \ref{lem4-2new}.
Since
$$
\int_{|y|<2}  V(y)w_i(y)dy\le C
\|V\|_{ L^{\frac n\alpha}(B_2) }
\|w_i\|_{ L^r(B_2) }
\le C
$$
for some $C$ independent of $i$,
we have
$$
\lim_{|z|\to 0}
\sup_i \|(L_iw_i)(\cdot+z)-(L_iw_i)(\cdot)\|_{ L^1(B_2) }=0.
$$
Therefore $\{L_iw_i\}$ is precompact
 in $L^1(B_2)$.

We know from Lemma \ref{lemfnew} that  $\{f_i\}$ 
converges to $f$ in $L^1(B_2)$.
So $\{w_i\}$ is  precompact  in $L^1(B_2)$.
After passing to a subsequence, $w_i\to w$
in $L^1(B_2)$. 
in view of  (\ref{qenew}),
$w\in L^r(B_2)$,  $w_i\to w$ in $L^p(B_2)$
for all $p<  r $,  $w^+\in  L^q(B_{\frac 12})$, and
\begin{equation}
\|w^+\|_{ L^q(B_{\frac 12}) }
\le  C\left(\|u\|_{ L^r(B_3) }+\|h\|_{ L^\nu(B_2) }\right).
\label{follow}
\end{equation}
  It follows that
$L_iw_i\to Lw$ in  $L^1(B_2)$.
Thus,
$$
w(x)= \int_{ |y|<2}
 \frac {
V(y)   w(y)
 }
{|x-y|^{n-\alpha}  }dy
+
f(x)+h(x)-\xi(x),  \qquad a.e.\ x\in B_2.
$$
Taking the difference of this  and (\ref{29-0new}),
we obtain
$$
(u-w)(x)
=
\int_{|y|<2}
 \frac {
V(y) (u-w)(y)  }
{|x-y|^{n-\alpha}  }dy,\quad a.e.\  x\in B_2.
$$
By the usual estimates and using $0<\delta(V)\le \bar\delta$
and $ C\bar \delta\le \frac 12$, 
$$
\|u-w\|_{
 L^r(B_2)
  }
\le C\bar \delta\|u-w\|_{
L^r(B_2)}
\le \frac 12
\|u-w\|_{  L^r(B_2)    }.
$$
It follows that
$
u=w$ a.e. in $ B_2$. 
Theorem \ref{thm5} follows from (\ref{follow}).

\vskip 5pt
\hfill $\Box$
\vskip 5pt

\noindent{\bf Proof of Corollary \ref{cor1}.}\ 
For $\epsilon>0$ small,
let
$$
u_\epsilon(x)=\epsilon^{ \frac{n-\alpha}2 }u(\epsilon x),
\ \ V_\epsilon(x)=\epsilon^\alpha V(\epsilon x),
\qquad x\in B_3,
$$
and
$$
h_\epsilon(x)=\epsilon^{ \frac{n-\alpha}2 }\int_{
3\epsilon <|y|< R_2 }
\frac {  V(y)u(y)  }
{ |\epsilon x-y|^{ n-\alpha}  }dy+ \epsilon^{ \frac{n-\alpha}2 }
h(\epsilon x).
$$
Then
$$
u_\epsilon(x)\le 
\int_{B_3} \frac {V_\epsilon(y) u_\epsilon(y) }{
|x-y|^{ n-\alpha }  }dy+h_\epsilon (x),
\qquad x\in B_2.
$$
Clearly, $u_\epsilon\in L^r(B_3)$ and
$h_\epsilon\in L^\nu(B_2)$.
Let $\bar\delta>0$ be the number in
Theorem \ref{thm5}, we fix some small  $\epsilon>0$ so that
$$
\|V_\epsilon\|_{ L^{\frac n\alpha}(B_3) }
=\|V\|_{ L^{\frac n\alpha}(B_{3\epsilon}) }<\bar\delta.
$$
Applying Theorem \ref{thm5} to $u_\epsilon$, we have
$u_\epsilon\in L^\nu(B_{\frac 12})$, i.e.
$u\in L^\nu(B_{\frac \epsilon2})$.

\vskip 5pt
\hfill $\Box$
\vskip 5pt

\noindent{\bf  Proof of Theorem \ref{thm1aa}.}\ 
Since $u\in L^{ \frac {2n}{n-\alpha} }_{loc}(\Bbb R^n)$, 
we have, by (\ref{1aa}), for some $|\bar x|<1$,
\begin{equation}
\int_{|y|> 2}
 \frac{ u(y)^{\frac {n+\alpha}{n-\alpha} }  }
{  |y|^{n-\alpha}  }dy
\le C
\int_{|y|> 2}
 \frac{ u(y)^{\frac {n+\alpha}{n-\alpha} }  }
{  |\bar x-y|^{n-\alpha}  }dy
\le \int_{\Bbb R^n} \frac{ u(y)^{\frac {n+\alpha}{n-\alpha} }  }
{  |\bar x-y|^{n-\alpha}  }dy
=u(\bar x)<\infty.
\label{23-1}
\end{equation}

For  any $R>0$,
we write
\begin{equation}
u(x)=I_R(x)+II_R(x):=\int_{|y|\le 2R}
 \frac{ u(y)^{\frac {n+\alpha}{n-\alpha} }  }
{  |x-y|^{n-\alpha}  }dy
+\int_{|y|> 2R}
 \frac{ u(y)^{\frac {n+\alpha}{n-\alpha} }  }
{  |x-y|^{n-\alpha}  }dy.
\label{Rnew}
\end{equation}
Take 
$$
V(x)=u(x)^{ \frac {2\alpha}{n-\alpha} },\qquad
h(x)= \int_{|y|> 2R}
 \frac{ u(y)^{\frac {n+\alpha}{n-\alpha} }  }
{  |x-y|^{n-\alpha}  }dy.
$$
Since
 $u\in L^{ \frac {2n}{n-\alpha} }_{loc}(\Bbb R^n)$,
$V\in L^{ \frac n\alpha}_{loc}(\Bbb R^n)$.
By (\ref{23-1}), $h\in L^\infty(B_R)$. 
For any $\nu>\frac n{n-\alpha}$, we have, by
Corollary \ref{cor1},
$u\in L^\nu( B_{ \epsilon(\nu) })$
for some $\epsilon(\nu)>0$.
Since any point can be taken as the origin, we have
proved that $u\in L^\nu_{loc}(\Bbb R^n)$
for all $1<\nu<\infty$.
By H\"older inequality, $I_R\in L^\infty(B_R)$. 
By (\ref{23-1}),  we can differentiate
$II_R(x)$ under the integral  for $|x|<R$, so $II_R\in C^\infty(\Bbb R^n)$.
 Since $R$ is
arbitrary,
$u\in L^\infty_{loc}(\Bbb R^n)$.
Back to (\ref{Rnew}),
$I_R$ is at least H\"older continuous
in $B_R$.  Since $R>0$ is arbitrary,
$u$ is  H\"older continuous
in $\Bbb R^n$.
Now  $u^{\frac {n+\alpha}{n-\alpha} }$ is H\"older continuous
in $
 B_{2R}$, the regularity of $I_R$
further improves and,  by bootstrap, we eventually have
$u\in C^\infty(\Bbb R^n)$.

\vskip 5pt
\hfill $\Box$
\vskip 5pt

 \section{ Proof of Theorem \ref{thm1new}}
In this section we prove Theorem \ref{thm1new}.
As shown in the last paragraph of Section 2,
$u\in C^\infty(\Bbb R^n)$.
By (\ref{1aa}) and the Fatou lemma,
\begin{equation}
\beta:=\liminf_{|x|\to \infty}(|x|^{n-\alpha}u(x) )
= \liminf_{|x|\to \infty}
 \int_{\Bbb R^n}  \frac{  |x|^{n-\alpha}u(y)^{\frac {n+\alpha}{n-\alpha} }  }
{ |x-y|^{ n-\alpha}  }dy
\ge \int_{\Bbb R^n}u(y)^{\frac {n+\alpha}{n-\alpha} }dy
>0.
\label{2prime11}
\end{equation}

For
  $x\in \Bbb R^n$,  $\lambda>0$, and  a positive function  $v$ on $\Bbb R^n$, 
let $
v_{x,\lambda}$ be given by (\ref{v24}).
Making a change of variables
$$
y=z^{x,\lambda}=x+  \frac{\lambda^2(z-x) }{  |z-x|^2 },
$$
we have
$$
dy= (\frac  \lambda {|z-x|})^{2n}
dz.
$$
Thus
\begin{eqnarray*}
\int_{ |y-x|\ge \lambda}
\frac {   v(y)^{\frac {n+\alpha}{n-\alpha} }  }
{  |\xi^{x,\lambda}-y|^{n-\alpha}
 }dy
&=& \int_{ |z-x|\le \lambda}
\frac {   v(z^{x,\lambda})^{\frac {n+\alpha}{n-\alpha} }  }
{  |\xi^{x,\lambda}-z^{x,\lambda}|^{n-\alpha}
 }(\frac  \lambda {|z-x|})^{2n}
dz\\
&=& \int_{ |z-x|\le \lambda}
\frac {   v_{x,\lambda}(z)^{\frac {n+\alpha}{n-\alpha} }  }
{  |\xi^{x,\lambda}-z^{x,\lambda}|^{n-\alpha}
 }(\frac  \lambda {|z-x|})^{n-\alpha}
dz.
\end{eqnarray*}
Since
$$
(\frac {|z-x|} \lambda )
(\frac {|\xi-x|}\lambda )
|\xi^{x,\lambda}-z^{x,\lambda}|=|\xi-z|,
$$
we have
\begin{equation}
(\frac \lambda {|\xi-x|})^{  {n-\alpha} }
\int_{ |y-x|\ge \lambda}
\frac {   v(y)^{\frac {n+\alpha}{n-\alpha} }  }
{  |\xi^{x,\lambda}-y|^{n-\alpha}
 }dy
= \int_{ |z-x|\le \lambda } \frac {  v_{x,\lambda}
(z)^{\frac {n+\alpha}{n-\alpha} }  }
{ |\xi-z|^{n-\alpha} } dz.
\label{811}
\end{equation}
Similarly,
\begin{equation}
(\frac \lambda {|\xi-x|})^{  {n-\alpha} }
\int_{ |y-x|\le \lambda}
\frac {   v(y)^{\frac {n+\alpha}{n-\alpha} }  }
{  |\xi^{x,\lambda}-y|^{n-\alpha}
 }dy
= \int_{ |z-x|\ge \lambda } \frac {  v_{x,\lambda}
(z)^{\frac {n+\alpha}{n-\alpha} }  }
{ |\xi-z|^{n-\alpha} } dz.
\label{911}
\end{equation}
For a positive solution $u$ of (\ref{1aa}),
applying 
(\ref{811}) and (\ref{911}) with $v=u$ and $v=u_{x,\lambda}$,
and using the fact $(\xi^{x,\lambda})^{x,\lambda}=\xi$ and
$(v_{x,\lambda})_{x,\lambda}\equiv v$, we obtain
\begin{equation}
u_{x,\lambda}(\xi)=
 \int_{ \Bbb R^n}  \frac {  u_{x,\lambda}
(z)^{\frac {n+\alpha}{n-\alpha} }  }
{ |\xi-z|^{n-\alpha} } dz, \qquad \forall\ \xi\in \Bbb R^n,
\label{1111}
\end{equation}
and
\begin{equation}
u(\xi)-u_{x,\lambda}(\xi)
=\int_{ |z-x|\ge \lambda}
K(x,\lambda; \xi,z)[ u(z)^{  \frac{n+\alpha}{n-\alpha} }-
u_{x,\lambda}(z)^{  \frac{n+\alpha}{n-\alpha} }]dz,
\label{1211}
\end{equation}
where
$$
K(x,\lambda; \xi,z)=\frac 1{ |\xi-z|^{n-\alpha} }
-(\frac \lambda{|\xi-x|} )^{n-\alpha}
 \frac 1{ |\xi^{x,\lambda}-z|^{n-\alpha} }.
$$
It is elementary to check that
$$
K(x,\lambda; \xi,z)>0,\qquad \forall\ |\xi-x|, |z-x|>\lambda>0.
$$

Formula (\ref{1111}) is the conformal invariance of the
integral equation (\ref{1aa}), see \cite{Lieb} and \cite{LL}.
\begin{lem}
For $x\in \Bbb R^n$, there exists $\lambda_0(x)>0$ such that
\begin{equation}
u_{x,\lambda}(y)\le u(y),\qquad
\forall\ 0<\lambda<\lambda_0(x), \ |y-x|\ge \lambda.
\label{1311}
\end{equation}
\label{lem111}
\end{lem}

\noindent{\bf Proof of Lemma \ref{lem111}.}\  The proof is essentially the same
as that of lemma 2.1 in \cite{LZhang}.  
For reader$'$s convenience, we include the
details.  Without loss of generality we may assume  $x=0$, 
and we use the notation $u_\lambda=u_{0,\lambda}$.

Since $\alpha<n$ and
$u$ is a positive $C^1$ function, there exists $r_0>0$ such that
$$
\nabla_y\left( |y|^{ \frac{n-\alpha}2 }u(y)\right)\cdot y>0,
\qquad\forall\ 0<|y|<r_0.
$$
Consequently
\begin{equation}
u_\lambda(y)<u(y), \qquad\forall\ 0<\lambda<|y|<r_0.
\label{1911}
\end{equation}
By (\ref{2prime11}) and the positivity and continuity of $u$,
\begin{equation}
u(z)\ge \frac 1{ C(r_0) |z|^{n-\alpha} }\qquad\forall\ |z|\ge r_0.
\label{1411}
\end{equation}
For small $\lambda_0\in (0, r_0)$ and for $0<\lambda<\lambda_0$,
$$
u_\lambda(y)
=(\frac \lambda{|y|})^{ n-\alpha}u(\frac {\lambda^2y}{|y|^2})
\le (\frac {\lambda_0}{|y|})^{ n-\alpha} \sup_{ B_{r_0}}u\le u(y),
 \qquad \forall\
|y|\ge r_0.
$$
Estimate (\ref{1311}), with $x=0$ and $\lambda_0(x)=\lambda_0$,
follows from (\ref{1911}) and the above.

\vskip 5pt
\hfill $\Box$
\vskip 5pt

Define, for $x\in \Bbb R^n$,
$$
\bar\lambda(x)=\sup\{\mu>0\ |\
u_{x,\lambda}(y)\le u(y)\ \forall\
0<\lambda<\mu, |y-x|\ge \lambda\}.
$$

\begin{lem} If $\bar\lambda(\bar x)<\infty$ for some
$\bar x\in \Bbb R^n$,  then
\begin{equation}
u_{\bar x, \bar\lambda(\bar x) }\equiv u\qquad
\mbox{on}\ \Bbb R^n.
\label{1711}
\end{equation}
\label{lem211}
\end{lem}

\noindent{\bf Proof of Lemma \ref{lem211}.}\
Without loss of generality, we may assume 
$\bar x=0$, and we use notations
$\bar \lambda=\bar\lambda(0), u_{\lambda}=
u_{0, \lambda}$.
By the definition of $\bar\lambda$,
\begin{equation}
u_{\bar \lambda}(y)\le u(y)\qquad\forall\ |y|
\ge \bar\lambda.
\label{bigger11}
\end{equation}
By (\ref{1211}), with $x=0$ and $\lambda=\bar\lambda$,
 and the positivity of the kernel,
either  $u_{\bar \lambda}(y)=u(y)$
for all $ |y|
\ge \bar\lambda$----then we are done-----or
$u_{\bar \lambda}(y)<u(y)$
 for all $ |y|> \bar\lambda$, which we assume below.
By
 the Fatou lemma, 
\begin{eqnarray*}
&&\liminf_{ |y|\to\infty}
|y|^{n-\alpha}(u-u_{\bar \lambda})(y)
\\
&=&
\liminf_{ |y|\to\infty}
\int_{|z|\ge \bar\lambda}
|y|^{n-\alpha}K(0, \bar\lambda;
y,z)[ u(z)^{ \frac{n+\alpha}{n-\alpha} }
-u_{\bar\lambda}(z)^{ \frac{n+\alpha}{n-\alpha} }]dz\\
&\ge &
 \int_{|z|\ge \bar\lambda}
\left( 1-(\frac {\bar\lambda}{ |z| })^{ n-\alpha }
\right)[ u(z)^{ \frac{n+\alpha}{n-\alpha} }
-u_{\bar\lambda}(z)^{ \frac{n+\alpha}{n-\alpha} }]dz>0.
\end{eqnarray*}
Consequently,
there exists $\epsilon_1\in (0,1)$ such that
$$
(u-u_{\bar \lambda})(y)\ge \frac { \epsilon_1}
{ |y|^{ n-\alpha} }\qquad
\forall\ |y|\ge \bar\lambda+ 1.
$$
By the above and the explicit formula of $u_\lambda$,
there exists $0<\epsilon_2<\epsilon_1$ such that
\begin{equation}
(u-u_{\lambda})(y)
\ge  \frac {\epsilon_1}
{ |y|^{ n-\alpha} }+(u_{\bar\lambda}-u_{\lambda})(y)
\ge  \frac {\epsilon_1}
{2 |y|^{ n-\alpha} }\ 
\forall\ |y|\ge \bar\lambda+ 1, 
\bar\lambda\le \lambda\le \bar \lambda+ \epsilon_2.
\label{1811}
\end{equation}

Now, for $\epsilon\in (0, \epsilon_2)$ which we choose below,
we have, for $\bar\lambda\le \lambda\le  \bar \lambda+ \epsilon$
and for $\lambda\le |y|\le  \bar\lambda+ 1$,
\begin{eqnarray*}
(u-u_\lambda)(y)
&=&
\int_{ |z|\ge \lambda}
K(0,\lambda; y,z)[ u(z)^{  \frac{n+\alpha}{n-\alpha} }-
u_{\lambda}(z)^{  \frac{n+\alpha}{n-\alpha} }]dz\\
&\ge& \int_{ \lambda\le |z|\le \bar\lambda+1}
K(0,\lambda; y,z)[ u(z)^{  \frac{n+\alpha}{n-\alpha} }-
u_{\lambda}(z)^{  \frac{n+\alpha}{n-\alpha} }]dz\\
&&+  \int_{  \bar\lambda+2\le |z|\le \bar\lambda+3}
K(0,\lambda; y,z)[ u(z)^{  \frac{n+\alpha}{n-\alpha} }-
u_{\lambda}(z)^{  \frac{n+\alpha}{n-\alpha} }]dz\\
&\ge& \int_{ \lambda\le |z|\le \bar\lambda+1}
K(0,\lambda; y,z)[ u_{\bar\lambda}(z)^{  \frac{n+\alpha}{n-\alpha} }-
u_{\lambda}(z)^{  \frac{n+\alpha}{n-\alpha} }]dz
\\
&&+  \int_{ \bar\lambda+2\le |z|\le \bar\lambda+3}
K(0,\lambda; y,z)[ u(z)^{  \frac{n+\alpha}{n-\alpha} }-
u_{\lambda}(z)^{  \frac{n+\alpha}{n-\alpha} }]dz.
\end{eqnarray*}

Because of (\ref{1811}), there exists $\delta_1>0$
such that
$$
u(z)^{  \frac{n+\alpha}{n-\alpha} }-
u_{\lambda}(z)^{  \frac{n+\alpha}{n-\alpha} }\ge \delta_1,
\quad \bar\lambda+2\le |z|\le  
\bar\lambda+3.
$$

Since
$$
 K(0,\lambda; y,z)=0,\qquad\forall\ |y|=\lambda,
$$ 
$$
\nabla_y K(0,\lambda; y,z)\cdot y\bigg|_{ |y|=\lambda }
=(n-\alpha)|y-z|^{ \alpha-n-2}(|z|^2-|y|^2)>0,
\quad \forall\  \bar\lambda+2\le |z|\le \bar\lambda+3,
$$
and  the function is smooth in the relevant region, 
we have, using also the positivity of the kernel,
$$
K(0,\lambda;y,z)\ge \delta_2(|y|-\lambda),
\forall\ \bar\lambda\le \lambda\le |y|\le \bar\lambda+1,
\bar\lambda+2\le |z|\le 
\bar\lambda+3,
$$
where $\delta_2>0$ is some constant independent of $\epsilon$.
It is easy to see that for some constant $C>0$  independent of $\epsilon$,
we have, for $\bar\lambda\le \lambda\le
\bar\lambda+\epsilon$, 
$$
| u_{\bar \lambda}(z)^{  \frac{n+\alpha}{n-\alpha} }-
u_{\lambda}(z)^{  \frac{n+\alpha}{n-\alpha} }|\le C(\lambda-\bar\lambda)
\le C\epsilon,
\quad\forall\  \bar\lambda\le \lambda\le |z|\le
\bar\lambda+1,
$$
and (recall that $\lambda\le |y|\le \bar\lambda+1$)
\begin{eqnarray*}
 \int_{ \lambda\le |z|\le \bar\lambda+1}
K(0,\lambda; y,z)dz&\le& 
|\int_{ \lambda\le |z|\le \bar\lambda+1}
\left(\frac 1{ |y-z|^{n-\alpha}  }-
\frac 1{ |y^\lambda-z|^{n-\alpha}  }
\right)dz|\\
&&+
\int_{ \lambda\le |z|\le \bar\lambda+1}
\bigg| (\frac \lambda{|y|})^{n-\alpha}-1\bigg|
\frac 1{ |y^\lambda-z|^{n-\alpha}  }dz
\\
&\le& C|y^\lambda-y|
+C(|y|-\lambda)\le  
C(|y|-\lambda).
\end{eqnarray*}

It follows from the above
that for  small $\epsilon>0$
we have, for $\bar\lambda\le \lambda\le \bar\lambda+\epsilon$ and
$\lambda\le |y|\le \bar\lambda+1$, 
\begin{eqnarray*}
(u-u_\lambda)(y)
&\ge &
-C\epsilon\int_{ \lambda\le |z|\le \bar\lambda+1}
K(0,\lambda; y,z)dz+ 
\delta_1 \delta_2(|y|-\lambda) \int_{ \bar \lambda+ 2\le |z|\le \bar\lambda+3}
dz\\
&\ge &
(\delta_1 \delta_2
 \int_{ \bar \lambda+ 2\le |z|\le \bar\lambda+3}
dz-C\epsilon)(|y|-\lambda)\ge 0.
\end{eqnarray*}
This and (\ref{1811}) violate the definition of
$\bar\lambda$. 
Lemma \ref{lem211} is established.

\vskip 5pt
\hfill $\Box$
\vskip 5pt

By the definition of $\bar \lambda(x)$,
$$
u_{x,\lambda}(y)\le u(y),\qquad\forall\
0<\lambda<\bar\lambda(x),  |y-x|\ge \lambda.
$$
Multiplying the above by $|y|^{n-\alpha}$ and sending
$|y|$ to infinity yields
\begin{equation}
\beta=\liminf_{|y|\to\infty}
|y|^{n-\alpha}u(y)\ge \lambda^{n-\alpha}u(x),
\qquad\forall\
0<\lambda<\bar\lambda(x).
\label{star11}
\end{equation}
On the other hand, if $\bar \lambda(\bar x)<\infty$, 
we use Lemma \ref{lem211} and multiply
  (\ref{1711}) by $|y|^{n-\alpha}$ and then send
$|y|$ to infinity to obtain
\begin{equation}
\beta=\lim_{|y|\to\infty}|y|^{n-\alpha}u(y)=
\bar\lambda(\bar x)^{n-\alpha}u(\bar x)
<\infty.
\label{starstar}
\end{equation}

\noindent{\bf Proof of Theorem \ref{thm1new}.}\
(i)\ If there exists some $\bar x\in \Bbb R^n$ 
such that $\bar\lambda(\bar x)<\infty$,
then, by (\ref{starstar}) and (\ref{star11}),
$\bar\lambda(x)<\infty$ for all $x\in \Bbb R^n$.
Applying Lemma \ref{lem211}, we have
$$
u_{x, \bar\lambda(x)}\equiv u\qquad\mbox{on}\ \Bbb R^n, \ \forall\ x\in
\Bbb R^n.
$$
By a calculus lemma (lemma 11.1 in \cite{LZhang}, see also lemma 2.5 
in \cite{LZhu} for $\alpha=2$), any $C^1$ positive function $u$ satisfying
the above must be of the form
(\ref{formnew}).

(ii) If $\lambda(x)=\infty$ for all $x\in \Bbb R^n$, then
$$
u_{x, \lambda}(y)\le u(y)\qquad
\forall\ |y-x|\ge \lambda>0, x\in \Bbb R^n.
$$
By another 
calculus lemma (lemma 11.2 in \cite{LZhang}, see also lemma 2.2 in 
\cite{LZhu} for $\alpha=2$), 
$u\equiv constant$, violating (\ref{1aa}).
Theorem \ref{thm1new} is established.

\vskip 5pt
\hfill $\Box$
\vskip 5pt

\section{Proof of Theorem \ref{thm1a}}

In this section we establish  Theorem \ref{thm1a}.
\begin{lem} For $n\ge 1$, $0<\alpha<n$ and $\mu>0$, let
$u$ be a Lebesgue measurable positive  solution of
(\ref{1a}) which is not identically equal to $\infty$. Then,
for any $t<\frac n{n-\alpha}$,
 $u\in L^\mu_{loc}(\Bbb R^n)\cap L^t_{loc}(\Bbb R^n)$,
\begin{equation}
\beta:=\liminf_{|x|\to \infty}(|x|^{n-\alpha}u(x) )
\ge \int_{\Bbb R^n}u(y)^\mu dy
>0,
\label{2prime}
\end{equation}
and
\begin{equation}
\int_{|y|>2} \frac {u(y)^\mu}{ |y|^{n-\alpha} } dy<\infty
\label{23}
\end{equation}
\label{lem3-1}
\end{lem}

\noindent{\bf Proof of Lemma \ref{lem3-1}.}\
Multiplying (\ref{1a}) by $|x|^{n-\alpha}$, we obtain
(\ref{2prime}) by applying the Fatou lemma.
Since $u$ is not identically equal to $\infty$, we see from 
(\ref{1a}) that
$u$ is finite almost everywhere.
So, for some  $x_1, x_2\in B_1$, $x_1\neq x_2$, we have
$$
\sum_{i=1}^2 \int_{ \Bbb R^n } 
\frac {  u(y)^\mu }{  |x_i-y|^{n-\alpha }  }dy
\le u(x_1)+u(x_2)<\infty.
$$ 
It follows that $u\in L^\mu_{loc}(\Bbb R^n)$ and (\ref{23}) holds.
For $R>0$, 
we write
\begin{equation}
u(x)=I_R(x)+II_R(x):=\int_{ |y|<2R} \frac { u(y)^\mu }
{ |x-y|^{n-\alpha} }dy+
\int_{ |y|>2R} \frac { u(y)^\mu }
{ |x-y|^{n-\alpha} }dy.
\label{R}
\end{equation}
Since $u\in L^\mu_{loc}(\Bbb R^n)$ and (\ref{23}) holds,
 $II_R\in L^\infty(B_R)$.
On the other hand, for any $1<t<\frac n{n-\alpha}$,
we have,  by Cauchy-Schwartz inequality,
\begin{eqnarray*}
\|I_R\|_{ L^t(B_R) }&\le&
 \int_{ |y|<2R }  u(y)^\mu \| |\cdot-y|^{\alpha-n}\|_{ L^t(B_R) }
dy\\
&\le&
 \| |\cdot-y|^{\alpha-n}\|_{ L^t(B_{3R})  }
 \int_{ |y|<2R } u(y)^\mu dy<\infty.
\end{eqnarray*}
Since $R>0$ is arbitrary, $u\in L^t_{loc}(\Bbb R^n)$.

\vskip 5pt
\hfill $\Box$
\vskip 5pt

\begin{lem} Assume   $n\ge 1$ and $0<\alpha<n$.

(i) For $0<\mu <\frac n{n-\alpha}$. Let  $u$ be 
a  positive
Lebesgue measurable solution of (\ref{1a}) which is not identically infinity.
Then $u\in C^\infty(\Bbb R^n)$.

(ii) For  $\mu \ge \frac n{n-\alpha}$.
 Let  $u\in  L^{  \frac {n (\mu-1) }\alpha  }_{loc}(\Bbb R^n)$ be 
a  positive
of (\ref{1a}).
Then $u\in C^\infty(\Bbb R^n)$.
\label{lem0a}
\end{lem}

\noindent{\bf Proof of Lemma \ref{lem0a}.}\
 
(i)\ For $0<\mu<\frac n{n-2}$. We know from Lemma \ref{lem3-1} that
$u\in L^t_{loc}(\Bbb R^n)$ for all
$t<\frac n{n-\alpha}$.
For any $R>0$, write $u$ as in (\ref{R}).
As usual, $II_R\in C^\infty(B_R)$.
For any 
$1<p<\frac n{ \mu(n-\alpha) }$,
let $\frac 1q=\frac 1p-\frac \alpha n$.
Then $q>\frac n{n-\alpha}$. By the property of
the Riesz potential,
$$
\|I_R\|_{ L^q(B_R) }
\le  
C\|u^\mu\|_{ L^p(B_{2R}) }=C\|u\|_{ L^{p\mu}(B_{2R}) }^\mu<\infty.
$$
So $u\in L^q_{loc}(\Bbb R^n)$.
Let $\mu'=\max(1, \mu)$,
since $u^\mu\le C+Cu^{\mu'}$,
we have
$$
u(x)\le C\int_{|y|<2R} \frac { V(y) u(y) }
{ |x-y|^{n-\alpha} }dy+h(x),\qquad x\in B_R,
$$
where
$$
V(y)= u(y)^{\mu'-1},\qquad h(x)=C+\int_{ |y|>2R} \frac { u(y)^\mu }
{ |x-y|^{n-\alpha} }dy.
$$
By (\ref{23}), $h\in L^\infty(B_R)$. Since $\frac {n(\mu'-1) }\alpha
<\frac n{n-\alpha}$,
$V\in L^{ \frac n\alpha }_{loc}(\Bbb R^n)$.
Since $u\in L^q_{loc}(\Bbb R^n)$ with $q>\frac n{n-\alpha}$,
 we have,
 by applying Corollary \ref{cor1},
$u\in L^\nu(B_{\epsilon(\nu)})$
for any $\nu>0$, where $\epsilon(\nu)>0$.  
Now, back to (\ref{R}), $I_R$ is
$C^\infty$ near the origin by bootstrapping.
By the translation invariance of the problem,
 $u\in C^\infty(\Bbb R^n)$.

(ii)\ For $\mu\ge \frac n{n-\alpha}$, 
let $V(y)=u(y)^{ \mu-1}$.  We know from Lemma \ref{lem3-1} that
$u\in L^t_{loc}(\Bbb R^n)$
for all $t<\frac n{n-\alpha}$.
Since $u\in L^{  \frac  {n(\mu-1)} \alpha  }_{loc}(\Bbb R^n)$
by the assumption, we also have $V\in L^{ \frac n\alpha }_{loc}(\Bbb R^n)$.
Now, for any $R_2>R_1>0$, let 
$$
h(y)=\int_{ |y|>R_2 }
\frac {  u(y)^\mu }
{  |x-y|^{ n-\mu} }dy.
$$
Then $u\in L^r(B_{R_2})$ with $r=\frac {n (\mu-1) }\alpha$,
 $V\in L^{ \frac n\alpha}(B_{R_2})$ 
$h\in L^\infty(B_{R_1})\subset L^\nu(B_{R_1})$ for any $\nu>r$,
and 
$$
u(x)=\int_{ |y|>R_2} \frac{  V(y) u(y) }{  |x-y|^{n-\alpha}  }dy
+h(x),\qquad x\in B_{R_1}.
$$
By Corollary \ref{cor1},
$u\in L^r(B_{R_1})$. 
Since $R_1>0$ is arbitrary, $u\in L^r_{loc}(\Bbb R^n)$
for all $r>1$.  Bootstrap as usual, $u\in C^\infty(\Bbb R^n)$.

\vskip 5pt
\hfill $\Box$
\vskip 5pt

For $x\in \Bbb R^n$,  $\lambda>0$
and a positive function $v$ on $\Bbb R^n$,
let $
v_{x,\lambda}$ be as in (\ref{v24}).

\begin{lem} For $n\ge 1$, $0<\alpha<n$ and $\mu>0$, 
let $u$ be a
positive solution of (\ref{1a}).  Then
\begin{equation}
u_{x,\lambda}(\xi)=
 \int_{ \Bbb R^n}  \frac {  u_{x,\lambda}
(z)^{\mu}  }
{ |\xi-z|^{n-\alpha} } 
\left( \frac\lambda{ |z-x|}\right)^{ n+\alpha-\mu(n-\alpha)}
dz, \qquad \forall\ \xi\in \Bbb R^n,
\label{11a}
\end{equation}
and
\begin{equation}
u(\xi)-u_{x,\lambda}(\xi)
=\int_{ |z-x|\ge \lambda}
K(x,\lambda; \xi,z)[ u(z)^{  \mu }-
\left( \frac\lambda{ |z-x|}\right)^{ n+\alpha-\mu(n-\alpha)}
u_{x,\lambda}(z)^{  \mu  }]dz,
\label{12a}
\end{equation}
where
$$
K(x,\lambda; \xi,z)=\frac 1{ |\xi-z|^{n-\alpha} }
-(\frac \lambda{|\xi-x|} )^{n-\alpha}
 \frac 1{ |\xi^{x,\lambda}-z|^{n-\alpha} }.
$$
Moreover,
$$
K(x,\lambda; \xi,z)>0,\qquad \forall\ |\xi-x|, |z-x|>\lambda>0.
$$
\label{lem3a}
\end{lem}

\noindent{\bf Proof of Lemma \ref{lem3a}.}\
The lemma for $\mu = \frac {n+\alpha}{n-\alpha}$
is established in Section 3.  The proof works for all
$\mu>0$ with minor modification.

\vskip 5pt
\hfill $\Box$
\vskip 5pt

\begin{lem}  For $n\ge 1$, $0<\alpha<n$ and $\mu>0$,
let $u\in C^1(\Bbb R^n)$ be a positive solution of
(\ref{1a}). Then for any  $x\in \Bbb R^n$, there exists $\lambda_0(x)>0$ such that
\begin{equation}
u_{x,\lambda}(y)\le u(y),\qquad
\forall\ 0<\lambda<\lambda_0(x), \ |y-x|\ge \lambda.
\label{13a}
\end{equation}
\label{lem1a}
\end{lem}

\noindent{\bf Proof of Lemma \ref{lem1a}.}\
This has been proved in Section 3 for 
$\mu=\frac {n+\alpha}{n-\alpha}$.
The same proof applies for all
$\mu>0$.

\vskip 5pt
\hfill $\Box$
\vskip 5pt

Define, for $x\in \Bbb R^n$,
$$
\bar\lambda(x)=\sup\{\mu'>0\ |\
u_{x,\lambda}(y)\le u(y)\ \forall\
0<\lambda<\mu', |y-x|\ge \lambda\}.
$$

\begin{lem} For $n\ge 1$, $0<\alpha<n$ and $0<\mu<\frac {n+\alpha}{n-\alpha}$,
let $u\in C^1(\Bbb R^n)$ be a positive solution of
(\ref{1a}).  Then
$\bar\lambda(x)=\infty$ for all $x\in \Bbb R^n$.
\label{lem5a}
\end{lem}

\noindent{\bf Proof of Lemma \ref{lem5a}.}\
We prove it by contradiction argument.
Suppose that $\bar\lambda(\bar x)<\infty$ for
some $\bar x\in \Bbb R^n$.
Without loss of generality, we may assume
$\bar x=0$, and we use notations
$\bar \lambda=\bar\lambda(0), u_{\lambda}=
u_{0, \lambda}$.
By the definition of $\bar\lambda$,
\begin{equation}
u_{\bar \lambda}(y)\le u(y)\qquad\forall\ |y|
\ge \bar\lambda.
\label{biggera}
\end{equation}
Since $n+\alpha-\mu(n-\alpha)>0$, 
$(\frac{\bar\lambda}{ |z|})^{ n+\alpha-\mu(n-\alpha)}
<1$ for $|z|>\bar\lambda$.
So, 
by (\ref{biggera}) and (\ref{12a}) with $x=0$ and $\lambda=\bar\lambda$,
 and the positivity of the kernel,
we have, for $|y|> \bar\lambda$,
\begin{eqnarray*}
(u-u_{\bar \lambda})(y)
&=&  \int_{|z|\ge \bar\lambda}
K(0, \bar\lambda;
y,z)[ u(z)^\mu-\left( \frac\lambda{ |z|}\right)^{ n+\alpha-\mu(n-\alpha)}
u_{\bar \lambda}(z)^\mu]dz\\
&\ge &
 \int_{|z|\ge \bar\lambda}
K(0, \bar\lambda;
y,z)[1-\left( \frac\lambda{ |z|}\right)^{ n+\alpha-\mu(n-\alpha)}
] u_{\bar \lambda}(z)^\mu dz>0.
\end{eqnarray*}
Thus, by  the Fatou lemma and the above,
\begin{eqnarray*}
&&\liminf_{ |y|\to\infty}
|y|^{n-\alpha}(u-u_{\bar \lambda})(y)
\\
&\ge &
 \liminf_{ |y|\to\infty}
\int_{|z|\ge \bar\lambda}
|y|^{n-\alpha}K(0, \bar\lambda;
y,z)[  u(z)^\mu-u_{\bar \lambda}(z)^\mu]dz\\
&\ge & 
\int_{|z|\ge \bar\lambda}
\left( 1- (\frac {\bar \lambda}{ |z| })^{n-\alpha}\right)
[  u(z)^\mu-u_{\bar \lambda}(z)^\mu]dz
>0.
\end{eqnarray*}
Consequently,
there exists $\epsilon_1\in (0,1)$ such that
$$
(u-u_{\bar \lambda})(y)\ge \frac { \epsilon_1}
{ |y|^{ n-\alpha} }\qquad
\forall\ |y|\ge \bar\lambda+ 1.
$$
By the above and the explicit formula of $u_\lambda$,
there exists $0<\epsilon_2<\epsilon_1$ such that
\begin{equation}
(u-u_{\lambda})(y)
\ge  \frac {\epsilon_1}
{ |y|^{ n-\alpha} }+(u_{\bar\lambda}-u_{\lambda})(y)
\ge  \frac {\epsilon_1}
{2 |y|^{ n-\alpha} }\
\forall\ |y|\ge \bar\lambda+ 1,
\bar\lambda\le \lambda\le \bar \lambda+ \epsilon_2.
\label{18a}
\end{equation}
Now, using (\ref{biggera}) and
(\ref{18a}) as in Section 3,
 for $\epsilon\in (0, \epsilon_2)$ which we choose below,
we have, for $\bar\lambda\le \lambda\le  \bar \lambda+ \epsilon$
and for $\lambda\le |y|\le  \bar\lambda+ 1$,
\begin{eqnarray*}
(u-u_\lambda)(y)
&\ge & \int_{ \lambda\le |z|\le \bar\lambda+1}
K(0,\lambda; y,z)[ u_{\bar\lambda}(z)^\mu-u_{\lambda}(z)^\mu]dz\\
&& +  \int_{  \bar\lambda+2\le |z|\le \bar\lambda+3}
K(0,\lambda; y,z)[ u(z)^\mu -u_{\lambda}(z)^\mu]dz.
\end{eqnarray*}
Because of (\ref{18a}), there exists $\delta_1>0$
such that
$$
u(z)^\mu -u_{\lambda}(z)^\mu\ge \delta_1,
\quad \bar\lambda+2\le |z|\le
\bar\lambda+3.
$$
It was shown in Section 3 that
$$
K(0,\lambda;y,z)\ge \delta_2(|y|-\lambda),
\forall\ \bar\lambda\le \lambda\le |y|\le \bar\lambda+1,
\bar\lambda+2\le |z|\le
\bar\lambda+3,
$$
where $\delta_2>0$ is some constant independent of $\epsilon$.
It is easy to see that for some constant $C>0$  independent of $\epsilon$,
we have, for $\bar\lambda\le \lambda\le
\bar\lambda+\epsilon$,
$$
| u_{\bar \lambda}(z)^\mu-
u_{\lambda}(z)^\mu|\le C(\lambda-\bar\lambda)
\le C\epsilon,
\quad\forall\  \bar\lambda\le \lambda\le |z|\le
\bar\lambda+1,
$$
and (recall that $\lambda\le |y|\le \bar\lambda+1$),
as in Section 3,
$$
 \int_{ \lambda\le |z|\le \bar\lambda+1}
K(0,\lambda; y,z)dz\le
C(|y|-\lambda).
$$
It follows from the above
that for  small $\epsilon>0$
we have, for $\bar\lambda\le \lambda\le \bar\lambda+\epsilon$ and
$\lambda\le |y|\le \bar\lambda+1$,
\begin{eqnarray*}
(u-u_\lambda)(y)
&\ge &
-C\epsilon\int_{ \lambda\le |z|\le \bar\lambda+1}
K(0,\lambda; y,z)dz+
\delta_1 \delta_2(|y|-\lambda) \int_{ \bar \lambda+ 2\le |z|\le \bar\lambda+3}
dz\\
&\ge &
(\delta_1 \delta_2
 \int_{ \bar \lambda+ 2\le |z|\le \bar\lambda+3}
dz-C\epsilon)(|y|-\lambda)\ge 0.
\end{eqnarray*}
This and (\ref{18}) violate the definition of
$\bar\lambda$.
 Lemma \ref{lem5a} is established.

\vskip 5pt
\hfill $\Box$
\vskip 5pt

\noindent{\bf Proof of Theorem \ref{thm1a}.}\
According to Lemma \ref{lem5a}, $\bar
 \lambda(x)=\infty$ for all $x\in \Bbb R^n$,
i.e.,
$$
u_{x, \lambda}(y)\le u(y)\qquad
\forall\ |y-x|\ge \lambda>0, x\in \Bbb R^n.
$$
By a calculus lemma (lemma 11.2 in \cite{LZhang}, see also lemma 2.2 in
\cite{LZhu} for $\alpha=2$),
$u\equiv constant$, violating (\ref{1a}).
Theorem \ref{thm1a} is established.

\vskip 5pt
\hfill $\Box$
\vskip 5pt

\section{Proof of Theorem \ref{thm1}}
In this section we establish Theorem \ref{thm1}.

\begin{lem} For $n\ge 1$, $p, q>0$, let $u$
be a non-negative Lebesgue measurable function in $\Bbb R^n$ satisfying
(\ref{1}).  Then
\begin{equation}
\int_{\Bbb R^n} (1+|y|^p)u(y)^{-q}dy<\infty,
\label{extra}
\end{equation}
\begin{equation}
\gamma:=\lim_{|x|\to\infty} |x|^{-p}u(x)=\lim_{|x|\to\infty}
\int_{\Bbb R^n}  \frac {|x-y|^p} { |x|^p } u(y)^{-q}dy
=\int_{\Bbb R^n} u(y)^{-q}dy \in (0,\infty),
\label{g}
\end{equation}
and, for  some constant $C\ge 1$,
\begin{equation}
 \frac {1+|x|^p}C\le u(x)\le C(1+|x|^p), \qquad \forall\ x\in \Bbb R^n.
\label{2}
\end{equation}
\label{lem0}
\end{lem}

\noindent{\bf Proof of 
Lemma \ref{lem0}.}\  We see from
(\ref{1}) that  $u$ must be positive everywhere and
$$
|\{y\in \Bbb R^n\ |\ u(y)<\infty\}|>0,
$$
where $|\ \cdot\ |$ denotes the Lebesgue measure of the set.
So there exist $R>1$
 and some  measurable set $E$ such that
$$
E\subset \{y\ |\ u(y)<R\}\cap B_{R},
$$
and
$$
|E|\ge \frac 1{R}.
$$

By (\ref{1}),
\begin{eqnarray*}
u(x)&=&\int_{\Bbb R^n}|x-y|^pu(y)^{-q}dy\ge 
\int_{E} |x-y|^pu(y)^{-q}dy\\
&\ge& (R)^{-q} \int_E  |x-y|^pdy, \qquad \forall\ x\in \Bbb R^n.
\end{eqnarray*}
The first inequality in (\ref{2}) follows from the above.

For some $1\le |\bar x|\le 2$,
$$
\int_{\Bbb R^n}|\bar x-y|^pu(y)^{-q}dy=u(\bar x)<\infty.
$$
We deduce
(\ref{extra}) from the  first inequality in (\ref{2}) and the above.  

For $|x|\ge 1$, 
$$
| \frac {|x-y|^p} { |x|^p } u(y)^{-q}|\le (1+|y|^p)u(y)^{-q},
$$
so, in view of (\ref{extra}),
 (\ref{g}) follows from the Lebesgue dominated convergence theorem.
The second inequality in (\ref{2}) follows from  
(\ref{1}), (\ref{extra}) and (\ref{g}).

\vskip 5pt
\hfill $\Box$
\vskip 5pt

\begin{lem}  For $n\ge 1$, $p, q>0$, let $u$
be a non-negative Lebesgue measurable function in $\Bbb R^n$ satisfying
(\ref{1}).  Then $u\in C^\infty(\Bbb R^n)$.
\label{lem0prime}
\end{lem}

\noindent{\bf Proof of Lemma \ref{lem0prime}.}\
For  $R>0$, writing
(\ref{1}) as
$$
u(x)=I_R(x)+II_R(x):=\int_{|y|\le 2R}
|x-y|^p u(y)^{-q}dy
+\int_{|y|> 2R}|x-y|^p u(y)^{-q}dy.
$$

Because of (\ref{extra}),  we can differentiate 
$II_R(x)$ under the integral  for $|x|<R$, and
therefore $II_R\in C^\infty(B_R)$.
On the other hand, since $u^{-q }\in
L^\infty(B_{2R})$, clearly 
$I_R$ is at least H\"older continuous
in $B_R$.  Since $R>0$ is arbitrary, 
$u$ is  H\"older continuous
in $\Bbb R^n$.
Now  $u^{-q }$ is H\"older continuous
in $
 B_{2R}$, the regularity of $I_R$  
further improves and,  by bootstrap, we eventually have  
 $u\in C^\infty(\Bbb R^n)$.
Lemma \ref{lem0prime} is established.

\vskip 5pt
\hfill $\Box$
\vskip 5pt

Let $v$ be a positive function on $\Bbb R^n$. For
  $x\in \Bbb R^n$ and $\lambda>0$,
consider 
$$
v_{x,\lambda}(\xi)=(\frac {|\xi-x|}\lambda)^{p}
v(\xi^{x,\lambda}),\qquad\xi\in \Bbb R^n,
$$
where
$$
\xi^{x,\lambda}=x+ \frac {\lambda^2(\xi-x)}{|\xi-x|^2}.
$$
Note that notation $v_{x,\lambda}$ in this section is different from
that in Section 1-4.

Making a change of variables
$$
y=z^{x,\lambda}=x+  \frac{\lambda^2(z-x) }{  |z-x|^2 },
$$
we have
$$
dy= (\frac  \lambda {|z-x|})^{2n}
dz.
$$
Thus
\begin{eqnarray*}
\int_{ |y-x|\ge \lambda} |\xi^{x,\lambda}-y|^p
 v(y)^{-q}dy
&=& \int_{ |z-x|\le \lambda} |\xi^{x,\lambda}-z^{x,\lambda}|^p
 v(z^{x,\lambda})^{-q}(\frac  \lambda {|z-x|})^{2n}
dz\\
&=& \int_{ |z-x|\le \lambda} |\xi^{x,\lambda}-z^{x,\lambda}|^p
(\frac  \lambda {|z-x|})^{2n-pq}
 v_{x,\lambda}(z)^{-q}
dz.
\end{eqnarray*}
Since 
$$
(\frac {|z-x|} \lambda )
(\frac {|\xi-x|}\lambda )
|\xi^{x,\lambda}-z^{x,\lambda}|=|\xi-z|,
$$
we have
\begin{eqnarray}
(\frac \lambda {|\xi-x|})^{-p }
&&\int_{ |y-x|\ge \lambda} |\xi^{x,\lambda}-y|^p v(y)^{-q}dy
\nonumber\\
&=& \int_{ |z-x|\le \lambda }|\xi-z|^p 
(\frac \lambda{ |z-x|})^{2n-pq+p}
v_{x,\lambda}
(z)^{-q}dz.
\label{8}
\end{eqnarray}
Similarly,
\begin{eqnarray}
&&(\frac \lambda {|\xi-x|})^{-p }
\int_{ |y-x|\le \lambda}
|\xi^{x,\lambda}-y|^p v(y)^{-q}dy
\nonumber\\
&=& \int_{ |z-x|\ge  \lambda }|\xi-z|^p 
(\frac \lambda{ |z-x|})^{2n-pq+p}
v_{x,\lambda}
(z)^{-q}dz.
\label{9}
\end{eqnarray}

\begin{lem} 
Let  $u$ be a positive
 solution  of (\ref{1}).
Then
\begin{equation}
u_{x,\lambda}(\xi)=
 \int_{ \Bbb R^n} |\xi-z|^p
(\frac \lambda{ |z-x|})^{2n-pq+p}
 u_{x,\lambda}
(z)^{-q}dz, \qquad \forall\ \xi\in \Bbb R^n,
\label{11}
\end{equation}
and
\begin{equation}
u_{x,\lambda}(\xi)-u(\xi)
=\int_{ |z-x|\ge \lambda}
k(x,\lambda; \xi,z)[ u(z)^{ -q }-
(\frac \lambda{ |z-x|})^{2n-pq+p}
u_{x,\lambda}(z)^{-q }]dz,
\label{12}
\end{equation}
where
$$
k(x,\lambda; \xi,z)=(\frac {|\xi-x|} \lambda)^{p}
|\xi^{x,\lambda}-z|^p-|\xi-z|^p.
$$
Moreover
$$
k(x,\lambda; \xi,z)>0,\qquad \forall\ |\xi-x|, |z-x|>\lambda>0.
$$
\label{lema}
\end{lem}

\noindent{\bf Proof of Lemma \ref{lema}.}\
Since  $(\xi^{x,\lambda})^{x,\lambda}=\xi$ and
$(v_{x,\lambda})_{x,\lambda}\equiv v$,
identity (\ref{11}) follows from (\ref{1}) and 
(\ref{8}) and (\ref{9}) with $v=u$. 
Similarly, using also (\ref{11}), 
\begin{eqnarray*}
u(\xi)&=& \int_{ |z-x|\ge \lambda}
|\xi-z|^p u(z)^{-q}dz+
\int_{ |y-x|<\lambda} |\xi-y|^p u(y)^{-q}dy\\
&=&  \int_{ |z-x|\ge \lambda}
|\xi-z|^p u(z)^{-q}dz\\
&&
+ (\frac {|\xi-x|}\lambda)^{p}
\int_{ |z-x|\ge \lambda}
|\xi^{x,\lambda}-z|^p
(\frac \lambda{ |z-x|})^{2n-pq+p}
 u_{x,\lambda}(z)^{-q}dz,
\end{eqnarray*}
\begin{eqnarray*}
u_{x,\lambda}(\xi)&=&
\int_{\Bbb R^n} |\xi-z|^p
(\frac \lambda{ |z-x|})^{2n-pq+p}
u_{x,\lambda}(z)^{-q}dz\\
&=&  
\int_{ |z-x|\ge \lambda}|\xi-z|^p
(\frac \lambda{ |z-x|})^{2n-pq+p}
u_{x,\lambda}(z)^{-q}dz\\
&&
 +(\frac {|\xi-x|}\lambda)^{p}
\int_{ |z-x|\ge \lambda} |\xi^{x,\lambda}-z|^p u(z)^{-q}dz.
\end{eqnarray*}
Identity (\ref{12}) follows from the above.
The positivity of the kernel $k$ is elementary.

\vskip 5pt
\hfill $\Box$
\vskip 5pt

\begin{lem} For $n\ge 1$, $p, q>0$,
 let $u$
be a  solution of (\ref{1}).  Then 
for any $x\in \Bbb R^n$, there exists $\lambda_0(x)>0$ such that
\begin{equation}
u_{x,\lambda}(y)\ge u(y),\qquad
\forall\ 0<\lambda<\lambda_0(x), \ |y-x|\ge \lambda.
\label{13}
\end{equation}
\label{lem1}
\end{lem}

\noindent{\bf Proof of Lemma \ref{lem1}.}\  The proof is similar to that 
 of lemma 2.1 in \cite{LZhang} and Lemma \ref{lem111} in Section 3.  
 Without loss of generality we may assume  $x=0$, 
and we use the notation $u_\lambda=u_{0,\lambda}$.

Since $p>0$ and
$u$ is a positive $C^1$ function, there exists $r_0>0$ such that
$$
\nabla_y\left( |y|^{ -\frac p2 }u(y)\right)\cdot y<0,
\qquad\forall\ 0<|y|<r_0.
$$
Consequently
\begin{equation}
u_\lambda(y)>u(y), \qquad\forall\ 0<\lambda<|y|<r_0.
\label{19}
\end{equation}

By (\ref{2}),
\begin{equation}
u(z)\le  C(r_0) |z|^p\qquad\forall\ |z|\ge r_0.
\label{14}
\end{equation}
For small $\lambda_0\in (0, r_0)$ and for $0<\lambda<\lambda_0$,
we have, using (\ref{2}) and (\ref{19}), 
$$
u_\lambda(y)
=(\frac{|y|} \lambda)^{p}u(\frac {\lambda^2y}{|y|^2})
\ge (\frac{|y|} {\lambda_0})^{p} \inf_{ B_{r_0}}u\ge u(y),
 \qquad \forall\
|y|\ge r_0.
$$
Estimate (\ref{13}), with $x=0$ and $\lambda_0(x)=\lambda_0$,
follows from (\ref{19}) and the above.

\vskip 5pt
\hfill $\Box$
\vskip 5pt

Define, for $x\in \Bbb R^n$,
$$
\bar\lambda(x)=\sup\{\mu>0\ |\
u_{x,\lambda}(y)\ge u(y)\ \forall\
0<\lambda<\mu, |y-x|\ge \lambda\}.
$$

\begin{lem}  For $n\ge 1$, $p>0$ and $0<q\le 1+\frac {2n}p$,
 let $u$
be a solution of (\ref{1}). 
Then
$$
\bar\lambda(x)<\infty,\qquad\forall\ x\in \Bbb R^n,
$$
and
\begin{equation}
u_{x, \lambda(x) }\equiv u\qquad
\mbox{on}\ \Bbb R^n, \qquad \forall\ x\in \Bbb R^n.
\label{17}
\end{equation}
Consequently, $q=1+\frac {2n}p$.
\label{lem2}
\end{lem}

\noindent{\bf Proof of Lemma \ref{lem2}.}\
By the definition of $\bar \lambda(x)$,
$$
u_{x,\lambda}(y)\ge u(y),\qquad\forall\
0<\lambda<\bar\lambda(x),  |y-x|\ge \lambda.
$$
Multiplying the above by $|y|^{-p}$ and sending
$|y|$ to infinity yields, using (\ref{g}), 
\begin{equation}
0<\gamma=\lim_{|y|\to\infty}
|y|^{-p}u(y)\le \lambda^{-p}u(x),
\qquad\forall\
0<\lambda<\bar\lambda(x).
\label{star}
\end{equation}
Thus  $\bar \lambda(x)<\infty$ for all $x\in \Bbb R^n$.

Now we prove (\ref{17}).
Without loss of generality, we may assume 
$x=0$, and we use notations
$\bar \lambda=\bar\lambda(0), u_{\lambda}=
u_{0, \lambda}$ and $y^\lambda=y^{0,\lambda}$.
By the definition of $\bar\lambda$,
\begin{equation}
u_{\bar \lambda}(y)\ge u(y)\qquad\forall\ |y|
\ge \bar\lambda.
\label{bigger}
\end{equation}

Since $2n-pq+p\ge 0$,
$(\frac {\bar \lambda}{ |z|})^{2n-pq+p}\le 1$
for $|z|\ge \bar \lambda$.
So, by (\ref{bigger}), 
(\ref{12}), with $x=0$ and $\lambda=\bar\lambda$,
 and the positivity of the kernel, 
either  $u_{\bar \lambda}(y)=u(y)$
for all $ |y|
\ge \bar\lambda$----then we are done (using (\ref{12})
to see that $2n-pq+p=0$)-----or
$u_{\bar \lambda}(y)>u(y)$
 for all $ |y|> \bar\lambda$, which we assume below.

By (\ref{12}),  with $x=0$ and $\lambda=\bar\lambda$, and the Fatou lemma
\begin{eqnarray*}
&&\liminf_{ |y|\to\infty}
|y|^{-p}(u_{\bar \lambda}-u)(y)
\\
&=&
\liminf_{ |y|\to\infty}
\int_{|z|\ge \bar\lambda}
|y|^{-p}
k(0, \bar\lambda;
y,z)
[ u(z)^{ -q }-(\frac {\bar \lambda}{ |z|})^{2n-pq+p}
u_{\bar\lambda}(z)^{  -q }]dz
\\
&\ge &
 \int_{|z|\ge \bar\lambda}
\left(
(\frac {|z|}{\bar\lambda})^p-1\right)
[ u(z)^{ -q }-u_{\bar\lambda}(z)^{  -q }]dz>0.
\end{eqnarray*}
Consequently, using also the positivity of $(u_{\bar\lambda}-u)$, 
there exists $\epsilon_1\in (0,1)$ such that
$$
(u_{\bar \lambda}-u)(y)\ge  \epsilon_1
 |y|^{p} \qquad
\forall\ |y|\ge \bar\lambda+ 1.
$$
By the above and the explicit formula of $u_\lambda$,
there exists $0<\epsilon_2<\epsilon_1$ such that
\begin{equation}
(u_{\lambda}-u)(y)
\ge  \epsilon_1
 |y|^{ p }+(u_{\lambda}-u_{\bar\lambda})(y)
\ge  \frac {\epsilon_1}
{2} |y|^{p }, \ 
\forall\ |y|\ge \bar\lambda+ 1, 
\bar\lambda\le \lambda\le \bar \lambda+ \epsilon_2.
\label{18}
\end{equation}

Recall that $2n-pq+p\ge 0$ and therefore 
$(\frac {\lambda}{ |z|})^{2n-pq+p}\le 1$
for $|z|\ge  \lambda$.
For $\epsilon\in (0, \epsilon_2)$ which we choose below,
we have, for $\bar\lambda\le \lambda\le  \bar \lambda+ \epsilon$
and for $\lambda\le |y|\le  \bar\lambda+ 1$,
\begin{eqnarray*}
(u_\lambda-u)(y)
&\ge &
\int_{ |z|\ge \lambda}
k(0,\lambda; y,z)[ u(z)^{-q  }-
u_{\lambda}(z)^{ -q }]dz\\
&\ge& \int_{ \lambda\le |z|\le \bar\lambda+1}
k(0,\lambda; y,z)[ u(z)^{ -q }-
u_{\lambda}(z)^{ -q }]dz\\
&&+  \int_{  \bar\lambda+2\le |z|\le \bar\lambda+3}
k(0,\lambda; y,z)[ u(z)^{-q }-
u_{\lambda}(z)^{  -q }]dz\\
&\ge& \int_{ \lambda\le |z|\le \bar\lambda+1}
k(0,\lambda; y,z)[ u_{\bar\lambda}(z)^{ -q }-
u_{\lambda}(z)^{  -q }]dz
\\
&&+  \int_{ \bar\lambda+2\le |z|\le \bar\lambda+3}
k(0,\lambda; y,z)[ u(z)^{  -q }-
u_{\lambda}(z)^{  -q }]dz.
\end{eqnarray*}

Because of (\ref{18}), there exists $\delta_1>0$
such that
$$
u(z)^{  -q }-
u_{\lambda}(z)^{  -q }\ge \delta_1,
\quad \bar\lambda+2\le |z|\le  
\bar\lambda+3.
$$

Since
$$
 k(0,\lambda; y,z)=0,\qquad\forall\ |y|=\lambda,
$$ 
$$
\nabla_y k(0,\lambda; y,z)\cdot y\bigg|_{ |y|=\lambda }
=p|y-z|^{ p-2}(|z|^2-|y|^2)>0,
\quad \forall\  \bar\lambda+2\le |z|\le \bar\lambda+3,
$$
and  the function is smooth in the relevant region, 
we have, using also the positivity of the kernel,
$$
k(0,\lambda;y,z)\ge \delta_2(|y|-\lambda),
\forall\ \bar\lambda\le \lambda\le |y|\le \bar\lambda+1,
\bar\lambda+2\le |z|\le 
\bar\lambda+3,
$$
where $\delta_2>0$ is some constant independent of $\epsilon$.
It is easy to see that for some constant $C>0$  independent of $\epsilon$,
we have, for $\bar\lambda\le \lambda\le
\bar\lambda+\epsilon$, 
$$
| u_{\bar \lambda}(z)^{ -q }-
u_{\lambda}(z)^{  -q }|\le C(\lambda-\bar\lambda)
\le C\epsilon,
\quad\forall\  \bar\lambda\le \lambda\le |z|\le
\bar\lambda+1,
$$
and (recall that $\lambda\le |y|\le \bar\lambda+1$)
\begin{eqnarray*} 
\int_{ \lambda\le |z|\le \bar\lambda+1}
k(0,\lambda; y,z)dz
&\le&
C(|y|-\lambda)+\int_{\lambda\le |z|\le \bar\lambda+1}
\left( |y^\lambda-z|^p-|y-z|^p\right)dz\\
&\le&
C(|y|-\lambda)+
 C|y^\lambda-y|
\le  
C(|y|-\lambda).
\end{eqnarray*} 

It follows from the above
that for  small $\epsilon>0$
we have, for $\bar\lambda\le \lambda\le \bar\lambda+\epsilon$ and
$\lambda\le |y|\le \bar\lambda+1$, 
\begin{eqnarray*}
(u_\lambda-u)(y)
&\ge &
-C\epsilon\int_{ \lambda\le |z|\le \bar\lambda+1}
k(0,\lambda; y,z)dz+ 
\delta_1 \delta_2(|y|-\lambda) \int_{ \bar \lambda+ 2\le |z|\le \bar\lambda+3}
dz\\
&\ge &
(\delta_1 \delta_2
 \int_{ \bar \lambda+ 2\le |z|\le \bar\lambda+3}
dz-C\epsilon)(|y|-\lambda)\ge 0.
\end{eqnarray*}
This and (\ref{18}) violate the definition of
$\bar\lambda$. 
Lemma \ref{lem2} is established.

\vskip 5pt
\hfill $\Box$
\vskip 5pt

\noindent{\bf Proof of Theorem \ref{thm1}.}\
According to  Lemma \ref{lem2},  $q=1+\frac{2n}p$ and
$$
u_{x, \bar\lambda(x)}\equiv u\qquad\mbox{on}\ \Bbb R^n, \ \forall\ x\in
\Bbb R^n.
$$
By a calculus lemma (lemma 11.1 in \cite{LZhang}, see also lemma 2.5 
in \cite{LZhu} for $\alpha=2$), any $C^1$ positive function $u$ satisfying
the above must be of the form
(\ref{form}).

\vskip 5pt
\hfill $\Box$
\vskip 5pt

\noindent{\bf Acknowledgment.}\ 
We thank  Chen,  Li and  Ou
for sending us  the preprint \cite{CLO} and thank 
 Ou
for an interesting talk on the work at a symposium.
We are grateful to Lieb for the encouragement
and for pointing out the
need to study the (essentially) uniqueness
of solutions of (\ref{1aa}) beyond
the  $L^\infty_{loc}(\Bbb R^n)$ class and bringing to our
attention of solutions  in \cite{Lieb} which are not maximizers.
This has led us
to study the regularity issue
and to establish Theorem \ref{thm1aa} and Theorem \ref{thm5}.  
We also thank Brezis for pointing out the
relation between Theorem \ref{thm5} and Lemma A.1 in
\cite{BL}.
The work was partially supported  by
NSF Grant DMS-0100819.

\end{document}